\documentclass[12pt,a4paper]{article}  % `article' with A4 paper size option
\usepackage{amsmath}
\usepackage{amsthm}                 % theorem extensions
\usepackage{amssymb}
\usepackage{graphicx}
\usepackage{psfrag}
\usepackage{caption}

\theoremstyle{plain}
\newtheorem{theorem}{Theorem}

\newtheorem{lemma}[theorem]{Lemma}

\newtheorem{proposition}[theorem]{Proposition}

\theoremstyle{definition}
\newtheorem{remark}[theorem]{Remark}
\newtheorem*{definition}{Definition}

\newcommand{\CC}{{\mathbb{C}}}
\newcommand{\HH}{{\mathbb{H}}}
\newcommand{\PP}{{\mathbb{P}}}
\newcommand{\QQ}{{\mathbb{Q}}}

\newcommand{\RR}{{\mathbb{R}}}
\newcommand{\ZZ}{{\mathbb{Z}}}

\newcommand{\calO}{{\cal O}}

\newcommand{\calI}{{\cal I}}

\newcommand{\calE}{{\cal E}}

\newcommand{\calX}{{\cal X}}
\newcommand{\calY}{{\cal Y}}
\newcommand{\calD}{{\cal D}}

\newcommand{\calS}{{\cal S}}

\newcommand{\bE}{{\calE}}
\newcommand{\e}[2]{{e_{#1}^{#2}}}
\newcommand{\be}[2]{{e_{#1}^{#2}}}
\newcommand{\E}[2]{{E_{#1}^{#2}}}
\newcommand{\F}[2]{{F_{#1}^{#2}}}

\newcommand{\Hom}{\mathrm{Hom}}
\newcommand{\Ext}{\mathrm{Ext}}
\newcommand{\Coh}{\mathrm{Coh}}

\newcommand{\TTT}{\mathsf{T}}
\newcommand{\chern}{\mathrm{ch}}
\newcommand{\Dsgr}{D_{\mathrm{sg}}^{\mathrm{gr}}}

\newcommand{\genfib}{Y} % generic fibre in the negatively smoothable Fuchs case

\title{McKay correspondence for the Poincar\'e series of Kleinian and Fuchsian singularities}
\date{}
\author{Wolfgang Ebeling and David Ploog
\thanks{
AMS Math.\ Subject Classification: 13D40, 18E30, 32S25..
}
}

\begin{document}

% global definition of labels in all figures
\psfrag{e11}{$\scriptstyle \E{1}{1}$}
\psfrag{e21}{$\scriptstyle \E{1}{2}$}
\psfrag{er1}{$\scriptstyle \E{1}{r}$}
\psfrag{er2}{$\scriptstyle \E{2}{r}$}
\psfrag{e1a1}{$\scriptstyle \E{\alpha_1-1}{1}$}
\psfrag{e2a2}{$\scriptstyle \E{\alpha_2-1}{2}$}
\psfrag{erar}{$\scriptstyle \E{\alpha_r-1}{r}$}
\psfrag{e}{$\scriptstyle E$}
\psfrag{u}{$\scriptstyle E-u$}
\psfrag{w}{$\scriptstyle u-w$}

\maketitle

\begin{abstract} We give a conceptual proof that the Poincar\'e series of the coordinate 
algebra of a Kleinian singularity and of a Fuchsian singularity of genus 0 is the 
quotient of the characteristic polynomials of two Coxeter elements. These Coxeter
elements are interpreted geometrically, using triangulated categories and spherical
twist functors.
\end{abstract}

\section*{Introduction}
Kleinian and Fuchsian singularities are defined by certain group actions: A Kleinian 
singularity is the quotient of $\CC^2$ by a finite subgroup $G\subset {\rm SU}(2)$, a 
Fuchsian singularity is defined by the action of a Fuchsian group (of the first kind)
$\Gamma \subset {\rm PSL}(2, \RR)$  on the cotangent bundle 
$T^{-1}_\HH$ of the upper half plane $\HH$. The classical McKay correspondence establishes 
a bijection between the irreducible representations of a finite subgroup 
$G \subset {\rm SU}(2)$ and the vertices of the affine extension of a Coxeter-Dynkin 
diagram of type ADE. The corresponding usual Coxeter-Dynkin diagram is the dual graph of 
the exceptional divisor of the minimal resolution of the singularity $(\CC^2/G,0)$. In 
\cite{Dolgachev07} a possible generalisation of this correspondence to Fuchsian groups of 
genus 0 is discussed. Here we show that there is a relation between invariants attached to 
the groups on one side and the Coxeter-Dynkin diagrams on  the other side.

To a group action there is associated a Poincar\'e series: The $k$-th coefficient of this 
infinite series is the dimension of the space of $G$-invariant homogeneous polynomials in 
two variables of degree $k$ in the Kleinian case and of the space of $\Gamma$-automorphic 
forms of weight $k$ on the upper half plane in the Fuchsian case. It was observed that 
these Poincar\'e series  can be described by the characteristic polynomials of the Coxeter 
elements attached to certain Coxeter-Dynkin diagrams.

In the Kleinian case, it was shown in \cite{E02}  that  the Poincar\'e series is the 
quotient of the characteristic polynomial of the Coxeter element of the corresponding 
affine Coxeter-Dynkin diagram by the one of the usual Coxeter-Dynkin diagram. In fact, 
this was derived from the McKay correspondence.

In \cite{E03} it was shown that for a Fuchsian singularity with $g=0$ the Poincar\'e 
series is also equal to the quotient of the characteristic polynomials of certain Coxeter 
elements. There this equality was proved by simply calculating both sides of this equation 
and showing that they are equal. No explanation could be given. In \cite{LP2} H.~Lenzing 
and J.~A.~de la Pe\~{n}a mention that this result also follows from results in 
\cite{L1, L2, LP1} in the realm of the representation theory of certain algebras related 
with a weighted projective line.

The aim of this paper is to give a uniform and conceptual geometrical proof of both the 
results for Kleinian and Fuchsian singularities. It is inspired by G.~Gonzalez-Sprinberg's 
and J.-L.~Verdier's geometric construction of the McKay correspondence (\cite{GV}, see 
also \cite{IN}) and the corresponding results in \cite{Dolgachev07} and  relies heavily on 
the general discussion in \cite{L2}. Moreover, we give a categorical interpretation of the 
results.

\section{Kleinian and Fuchsian singularities} \label{KandF}
A Kleinian or Fuchsian singularity is a normal surface singularity $(X,x)$ with a good 
$\CC^\ast$-action. This means that $X={\rm Spec}(A)$ is a normal two-dimensional affine 
algebraic variety  over $\CC$ which is smooth outside its {\em vertex} $x$. Its coordinate 
ring $A$ has the structure of a graded $\CC$-algebra $A = \bigoplus_{k=0}^\infty A_k$, $A_0=\CC$, 
and $x$ is defined by the maximal ideal $\mathfrak{m}= \bigoplus_{k=1}^\infty A_k$. 

A natural compactification of $X$ is given by $\overline{X}:= {\rm Proj}(A[t])$, 
where $t$ has degree 1 for the grading of $A[t]$ (see \cite{Pinkham77a}). This is a 
normal projective surface with $\CC^*$-action, and $\overline{X}$ may acquire additional
singularities on the boundary $\overline{X}_\infty:=\overline{X}\setminus X={\rm Proj}(A)$.
Let $\pi:S\to\overline{X}$ be the minimal normal crossing resolution of all singularities 
of $\overline{X}$. Thus $S$ is a smooth, projective surface.

For a normal surface singularity $(X,x)$, we have the following characterisations:
$(X,x)$ is Kleinian if it is a rational double point.
And $(X,x)$ is Fuchsian if the canonical sheaf of $\overline{X}$
is trivial. In this case, the singularities on the boundary are all of type $A_\mu$.
The genus of the Fuchsian singularity is defined as the genus $g=g(\overline{X}_\infty)$ 
of the boundary.

\paragraph{Equivariant description of the singularities.}
%Following I.~Dolgachev \cite{Dolgachev75}, we give another description of these types of 
%singularities. 
Following I.~Dolgachev \cite{Dolgachev75}, there exist a simply connected Riemann surface
$\cal D$, a discrete cocompact subgroup $\Gamma$ of $\mbox{Aut}({\cal D})$ and a line 
bundle $\cal L$ on $\cal D$ to which the action of $\Gamma$ lifts such that 
 $A_k = H^0({\cal D}, {\cal L}^k)^\Gamma$.
Let $Z:={\cal D}/\Gamma$. By \cite[Theorem~5.1]{Pinkham77a} (see also
\cite[Theorem~5.4.1]{Wagreich81}), there exist a divisor
$D_0$ on $Z$, points $p_1, \ldots , p_r \in Z$, and integers $\alpha_i$, $\beta_i$
with $0<\beta_i < \alpha_i$ and $(\alpha_i,\beta_i)=1$ for 
$i=1, \ldots, r$ with
\[ A_k =L (D^{(k)}), \quad 
   D^{(k)} := kD_0 + \sum_{i=1}^r \left[ k \frac{\alpha_i - \beta_i}{\alpha_i} \right] p_i
  \quad\mbox{ for } k\geq0 . \]
Here $[x]$ denotes the largest integer $\leq x$, and $L(D):=H^0(Z,\calO_Z(D))$ for a 
divisor $D$ on $Z$ denotes the linear space of meromorphic functions $f$ on $Z$ such 
that $(f) \geq -D$. We enumerate the points $p_i$ so that
 $\alpha_1\leq \alpha_2 \leq \ldots \leq \alpha_r$. 
 The genus $g$ of $Z$ coincides with the genus of $\overline{X}_\infty$. We define 
$b :={\rm degree}\, D_0 +r$. Then 
 $\{ g; b; (\alpha_1, \beta_1), \ldots , (\alpha_r, \beta_r)\}$ 
are called the {\em orbit invariants} of $(X,x)$, cf.\ e.g.\ \cite{Wagreich83}. Define
$\mbox{vdeg}({\cal L}):= -b + \sum_{i=1}^r \frac{\beta_i}{\alpha_i}$.

Assume that $(X,x)$ is Gorenstein. By \cite{Dolgachev83}, there exists an integer
$R$ such that ${\cal L}^{-R}$ and the tangent bundle $T_{\cal D}$ of $\cal D$ are
isomorphic as $\Gamma$-bundles and
\begin{eqnarray*}
R \cdot \mbox{vdeg}({\cal L}) & = & 2 -2g - r + \sum_{i=1}^r \frac{1}{\alpha_i}, \\
R\beta_i & \equiv & 1 \ \mbox{mod} \, \alpha_i, \quad i=1, \ldots , r.
\end{eqnarray*}

If $R=-1$ then ${\cal D} = \PP^1(\CC)$, $g=0$, and $\Gamma$ is a finite subgroup of 
$\mbox{Aut}(\PP^1(\CC)) = PGL(2,\CC)$. In this case, $(X,x)$ is a Kleinian singularity. 
The orbit invariants are 
 $ \{ 0; 2; (\alpha_1, \alpha_1-1), \ldots , (\alpha_r, \alpha_r -1) \} $
and $\sum_{i=1}^r \frac{1}{\alpha_i} > r-2$ (cf.\ \cite[(3.2)]{Wagreich83}). This implies $r \leq 3$.
In this case
\[ D^{(k)} = D^{(k)}_{\rm Klein} := kD_0 + \sum_{i=1}^r \left[ \frac{k}{\alpha_i} \right] p_i, \]
where  $D_0$ is a divisor on $Z$ of degree $2-r$.

If $R=1$, then ${\cal D} = \HH$ and $\Gamma \subset PSL(2,\RR)$ is
a finitely generated cocompact {\em Fuchsian group of the first kind}. This means that
$\Gamma$ acts properly discontinuously on $\HH$ and that the quotient $Z =\HH / \Gamma$
is a compact Riemann surface. The divisor $D_0$ is the canonical divisor, the points 
$p_1, \ldots, p_r \in Z$ are the branch points of the map $\HH \to Z$,  $\alpha_i$ is 
the ramification index over
$p_i$, and $\beta_i=1$ for $i=1, \ldots ,r$. Hence the orbit invariants are
 $ \{ g; 2g-2+r; (\alpha_1,1), \ldots, (\alpha_r,1) \} $
and $\sum_{i=1}^r \frac{1}{\alpha_i} < r+2g-2$ (cf.\ \cite[(1.7)]{Looijenga84}). In the case $g=0$ this 
implies $r \geq 3$.
According to \cite{Looijenga84}, $(X,x)$  a Fuchsian singularity. Here
\[ D^{(k)} = D^{(k)}_{\rm Fuchs} := kD_0 + 
   \sum_{i=1}^r \left[ k \frac{\alpha_i - 1}{\alpha_i} \right] p_i \]
and ${\rm deg}(D_0)=2g-2$.

\paragraph{Configurations of exceptional curves.} Fuchsian and Kleinian singularities contain
two particular curve configurations, the roles of which are interchanged, however.

Let $(X,x)$ be a Kleinian singularity. Then the minimal resolution of the 
singularity $x$ has an exceptional divisor with the dual graph $\bE$ depicted 
in Figure~\ref{Fig1}. Here all vertices correspond to rational curves of 
self-intersection $-2$, the mutual intersection numbers are either 0 or 1, 
and two vertices are joined by an edge if and only if the intersection number 
of the corresponding rational curves is equal to 1. The surface $S$ also 
contains a system $\bE'$ of rational curves intersecting as in Figure~\ref{Fig2}.
This system is the boundary divisor of the compactification, i.e.\ the preimage
of $\overline{X}_\infty$ under $\pi:S\to\overline{X}$.

\vspace{2.5ex}

\noindent
\parbox[t]{0.45\textwidth}{
  \includegraphics[width=\linewidth]{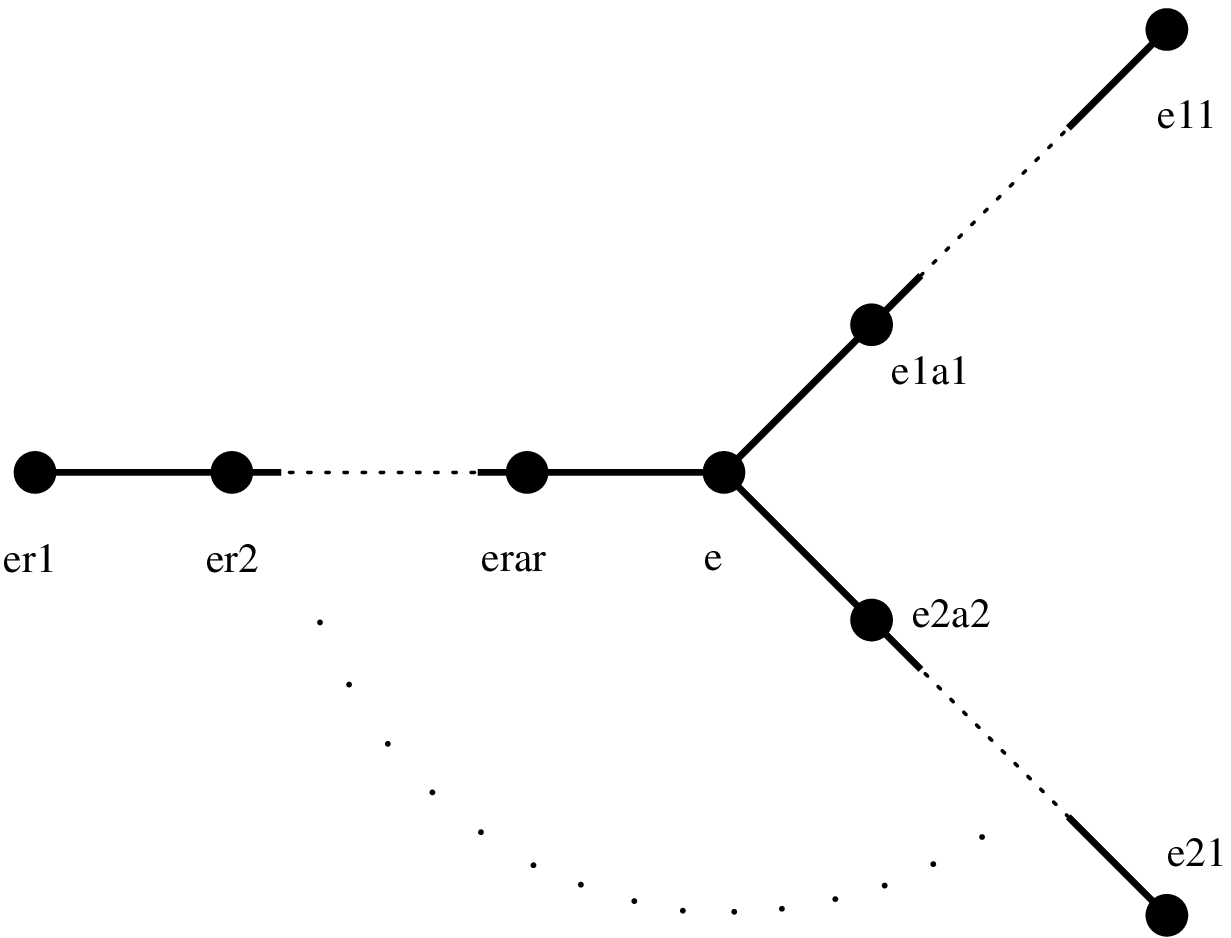}
  \captionof{figure}{The dual graph of the configuration $\bE$. Vertices 
                     correspond to rational curves of self-intersection $-2$;
                     the curves are labelled as in the figure.}
  \label{Fig1}
}
\hfill
\parbox[t]{0.45\textwidth}{
  \psfrag{a0}{$\scriptstyle-1$}
  \psfrag{a1}{$\scriptstyle-\alpha_1$}
  \psfrag{a2}{$\scriptstyle-\alpha_2$}
  \psfrag{a3}{$\scriptstyle-\alpha_3$}
  \psfrag{a4}{$\scriptstyle-\alpha_4$}
  \psfrag{a5}{$\scriptstyle-\alpha_5$}
  \psfrag{ar}{$\scriptstyle-\alpha_r$}
  \centering
  \includegraphics[width=0.75\linewidth]{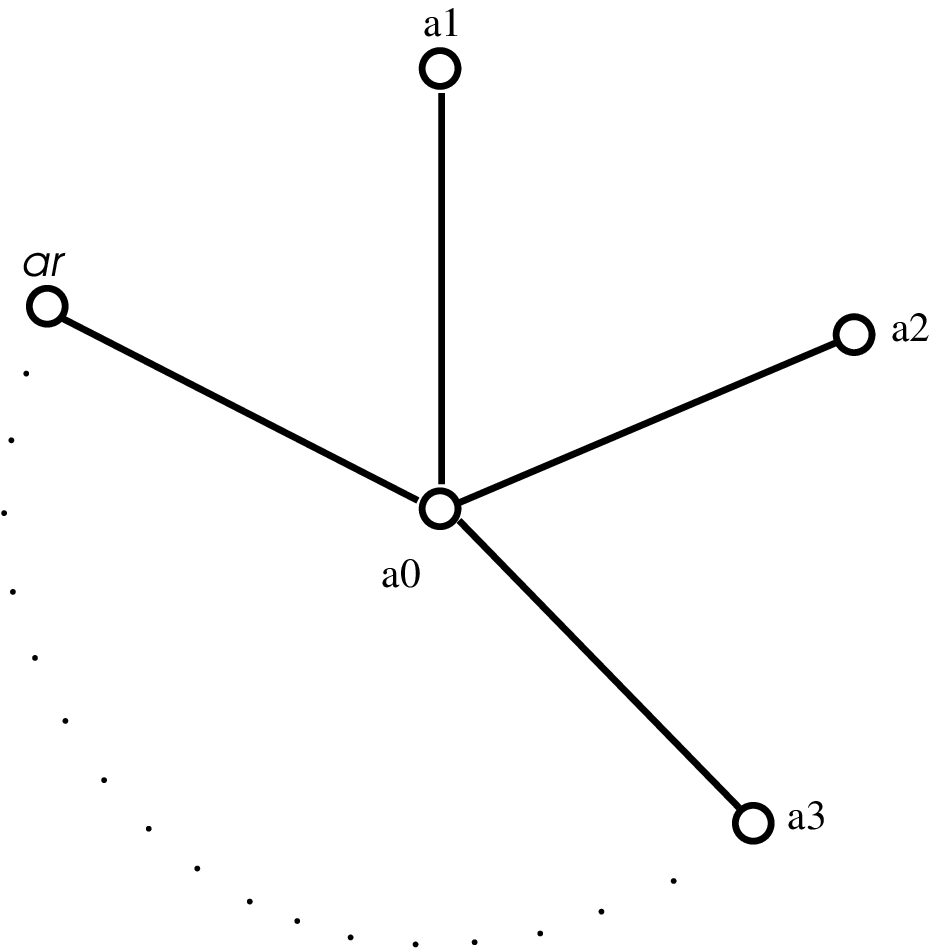}
  \captionof{figure}{The dual graph of the configuration $\bE'$. Each vertex 
                     corresponds to a curve of negative self-intersection 
                     indicated by the label.}
  \label{Fig2}
}

\vspace{1.5ex}

Now assume that $(X,x)$ is Fuchsian with $g=0$. In this case the minimal 
normal crossing resolution of $x$ looks as in Figure~\ref{Fig2}. The variety 
$\overline{X}$ has $r$ cyclic quotient singularities of type
$(\alpha_1,\alpha_1-1), \ldots , (\alpha_r,\alpha_r-1)$
along $\overline{X}_\infty:=\overline{X}-X$ \cite[Lemma~4.1]{Pinkham77a}. A cyclic
quotient singularity of type $(\alpha,\alpha-1)$ is a singular point of type
$A_{\alpha-1}$. 
The preimage $\widetilde{X}_\infty$ of $\overline{X}_\infty$ under $\pi : S \to \overline{X}$
consists of the strict transform $E$ of $\overline{X}_\infty$ and $r$ chains
$\E{1}{i}, \ldots , \E{\alpha_i-1}{i}$, $i=1, \ldots , r$, of rational curves of
self-intersection $-2$ which intersect again according to the dual graph shown in
Figure~\ref{Fig1}. By the adjunction formula and $g=0$, the self-intersection number of 
the curve $E$ is also $-2$.

%%%%%%%%%%%%%%%%%%%%%%%%%%%%%%%
\section{Lattices and the  Poincar\'e series of Kleinian and Fuchsian singularities}
We consider the {\em
Poincar\'{e} series} of the algebra $A$ 
\[ p_A(t)= \sum_{k=0}^\infty \dim(A_k) t^k . \]
We shall give a description of $p_A$ for the Kleinian singularities and the Fuchsian 
singularities with $g=0$; before that, we recall the notions of root and Coxeter element.

\paragraph{Coxeter elements and lattice extensions.}
If $(V,\langle-,-\rangle)$ is an arbitrary lattice and $e \in V$ is a {\em root},
i.e.\ $\langle e,e \rangle=-2$, then the reflection corresponding to $e$ is 
defined by 
\[ s_e(x) = x - \frac{2\langle x,e \rangle}{\langle e, e \rangle} e 
          = x + \langle x,e\rangle e \quad \mbox{for } x \in V. \]
If $B=(e_1, \ldots , e_n)$ is an ordered basis consisting of roots, then the 
{\em Coxeter element} $\tau$ corresponding to $B$ is defined by
\[ \tau= s_{e_1} s_{e_2} \cdots s_{e_n}. \] 
For a Coxeter element $\tau$, denote its characteristic polynomial by 
$\Delta(t)= \det (\tau-t \cdot {\rm id})$.

We consider the following two extensions of $V$. Let $U$ be a unimodular 
hyperbolic plane, i.e.\ a free $\ZZ$-module with basis $u,w$ and symmetric 
bilinear form defined by 
\[ \langle u,u \rangle= \langle w,w \rangle =0, \quad 
    \langle u,w \rangle =\langle w,u \rangle = 1. \]
Let $V[u]:=V \stackrel{\perp}{\oplus} \ZZ u$ and $V[u,w]:=V \stackrel{\perp}{\oplus} U$. 

\paragraph{Application to the geometric case.}
Let $V_-$ be the abstract lattice associated to Figure~\ref{Fig1}, i.e.\ the free 
$\ZZ$-module spanned by the ordered basis 
 $B_-:=(\E{1}{1}, \ldots, \E{\alpha_1-1}{1};  \ldots ; \E{1}{r}, \ldots , \E{\alpha_r-1}{r}; E)$ 
with bilinear form $\langle-,-\rangle$ given by the intersection numbers. 
We observe that all vectors of the basis square to $-2$; denote by 
 $\tau_-\in\mathrm{O}(V_-,\langle-,-\rangle)$
the Coxeter element of $V_-$ corresponding to the above basis.

\begin{figure}
\includegraphics[scale=0.5]{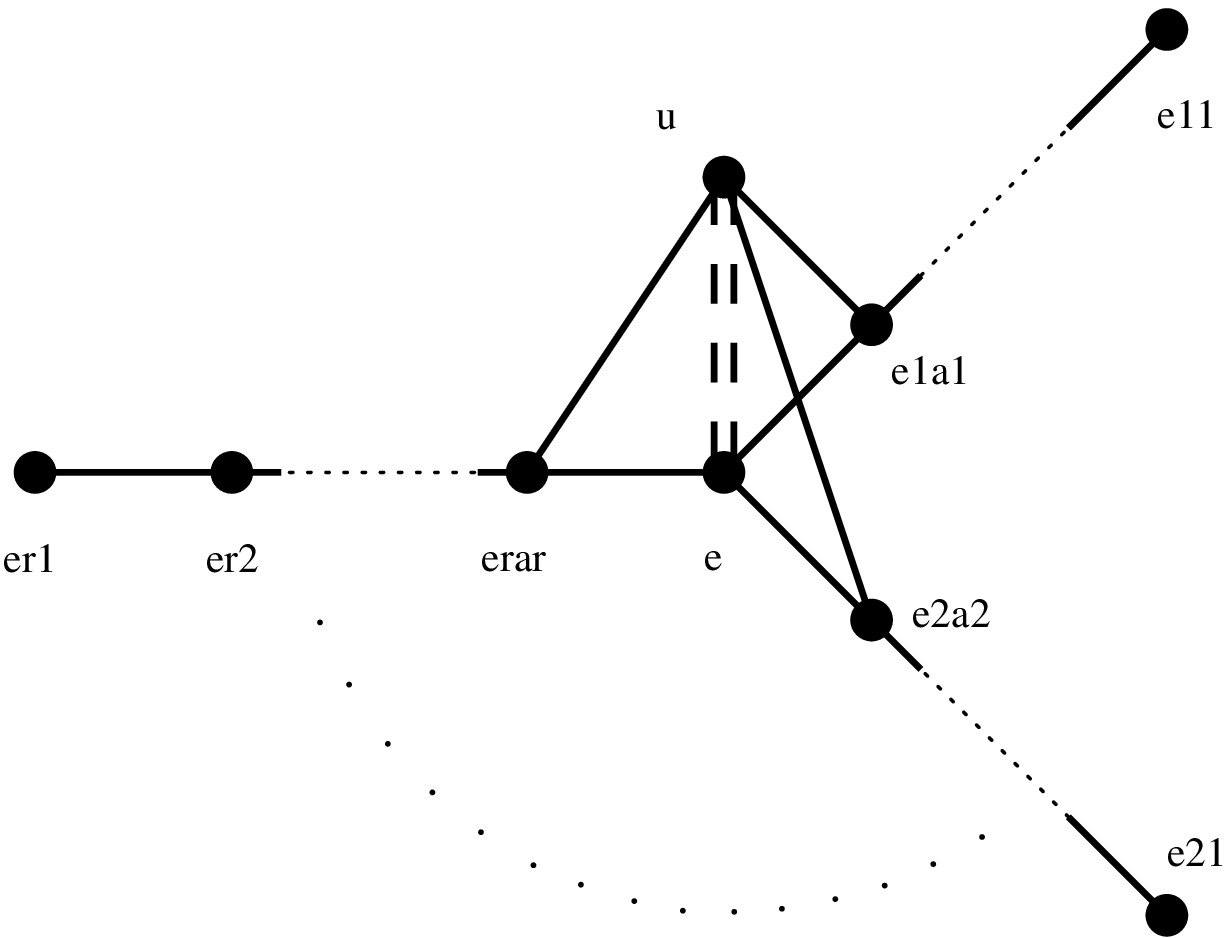}
\hfill
\includegraphics[scale=0.5]{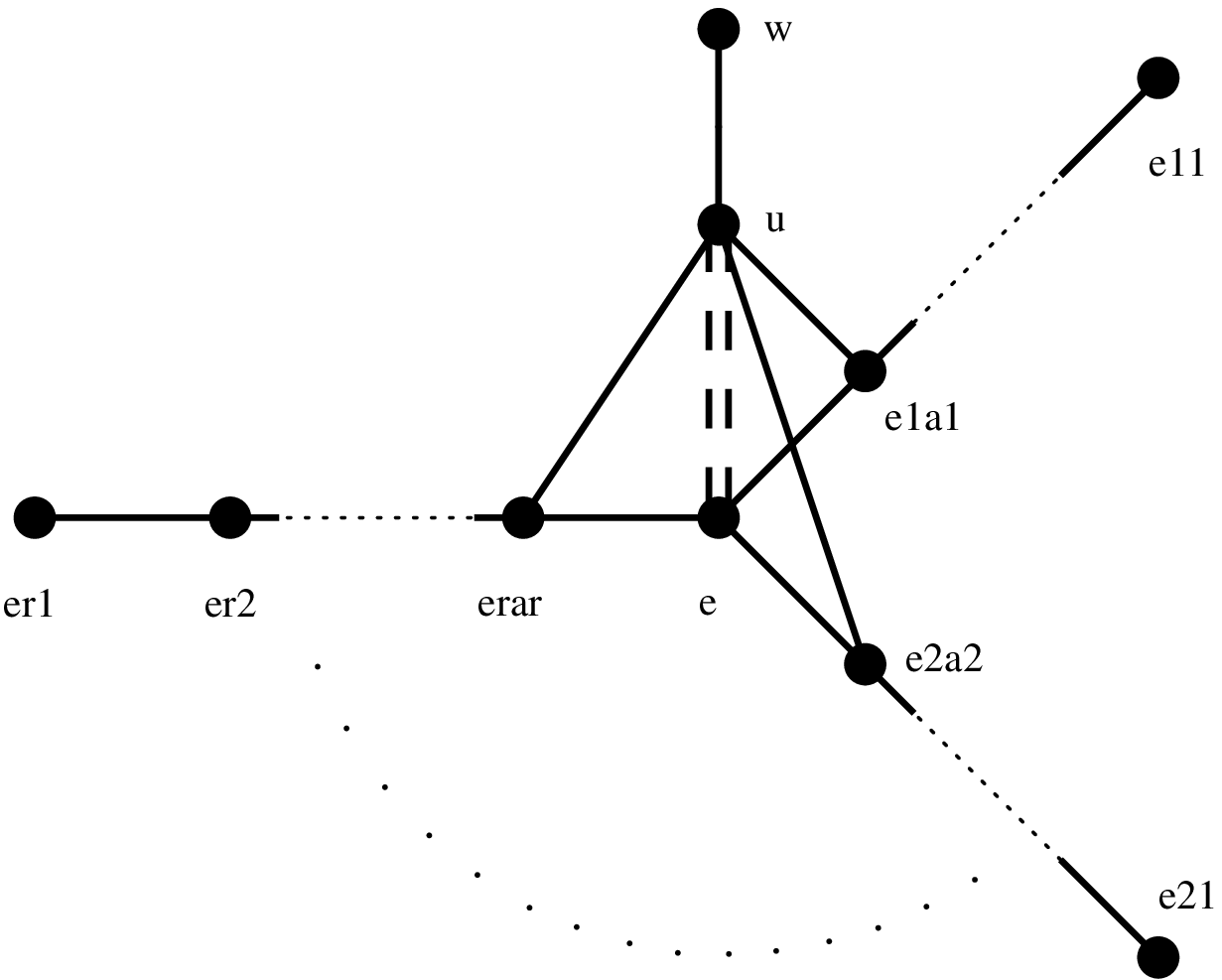}
\caption{The Coxeter-Dynkin diagrams of $V_0$ (left) and $V_+$ (right) } 
\label{Fig3}
\end{figure}

Let $V_0:=V_-[u]$ with basis $B_0:=(B_- , E-u)$ and $V_+:=V_-[u,w]$ with basis 
$B_+:=(B_-, E-u, u-w)$. See Figure~\ref{Fig3} for the graphs.
Let $\tau_-$, $\tau_0$, $\tau_+$ be the Coxeter elements 
corresponding to the bases $B_-$, $B_0$, $B_+$, and let $\Delta_-, \Delta_0, \Delta_+$ be 
their characteristic polynomials, respectively.

% $$ \xymatrix{ 
%  & & & & & & {\bullet} \ar@{}^{e_1^{(1)}}[dr] & \\
%  & & & & & {\cdots} \ar@{-}[ur]  & & \\
%  & & & {\bullet} \ar@{==}[d] \ar@{-}[r]  \ar@{-}[rdd] \ar@{}^{\widetilde{e}}[u] & {\bullet} \ar@{-}[ur] \ar@{}^{e_{\alpha_1-1}^{(1)}}[dr] & & & \\
%  {\bullet} \ar@{-}[r] \ar@{}_{e_1^{(r)}}[d]  & {\cdots} \ar@{-}[r]  & {\bullet} \ar@{-}[r] \ar@{-}[ur]   \ar@{}_{e_{\alpha_r-1}^{(r)}}[d] & {\bullet} \ar@{-}[ur] \ar@{-}[dr] \ar@{}_{e}[d] & & & & \\
% & & & & {\bullet} \ar@{-}[dr] \ar@{}_{e_{\alpha_2-1}^{(2)}}[ur] & & &  \\
%  & & & & & {\cdots} \ar@{-}[dr] & & \\
% & & & & & & {\bullet}  \ar@{}_{e_1^{(2)}}[ur] & 
%   } $$

\begin{theorem} \label{thm:main}
\begin{itemize}
\item[{\rm (i)}] For a Kleinian singularity we have
\[ p_A = \frac{\Delta_-}{\Delta_0}. \]
\item[{\rm (ii)}] For a Fuchsian singularity with $g=0$ we have
\[ p_A = \frac{\Delta_+}{\Delta_0}. \]
\end{itemize}
\end{theorem} 

Theorem~\ref{thm:main}~(i) was derived in \cite{E02} from the McKay correspondence
for the case $(X,x)$ not of type $A_{2n}$. If $(X,x)$ is a Kleinian singularity
of type $A_\mu$, then the weights and the degree are not unique. In \cite{E02} it
was tacitly assumed that $A_\mu$ is given by the standard equation
$f(x,y,z)=x^{\mu+1}+y^2+z^2=0$.

Theorem~\ref{thm:main}~(ii) follows from \cite[Proposition~1 and Remark~1]{E03}. Part (ii) 
was also proven independently by H.~Lenzing and J.~A.~de la Pe\~{n}a (see \cite{LP2} 
and the references there).

\begin{remark} The Fuchsian singularities with $g=0$, $r=3$, and where the
algebra $A$ is generated by 3 elements, are precisely the 14 exceptional
unimodal hypersurface singularities in V.~I.~Arnold's classification. In
this case there exists a distinguished basis of vanishing cycles of the
Milnor lattice such that the Coxeter-Dynkin diagram is equal to the graph 
of Figure~\ref{Fig2} for another triple of numbers $(p_1,p_2,p_3)$ (the 
{\em Gabrielov numbers}) which is in general different from the triple
$(\alpha_1,\alpha_2,\alpha_3)$ (which are called {\em Dolgachev numbers}).
Arnold observed that there is an involution on the set of these
singularities ({\em strange duality}) which interchanges Dolgachev and
Gabrielov numbers. This duality can be considered as part of the mirror
symmetry of K3 surfaces (for details see e.g.\ \cite{E03}). In this case, 
the Coxeter element $\tau_+$ corresponding to the triple $(p_1,p_2,p_3)$ 
is the monodromy operator of the singularity.
\end{remark}

%%%%%%%%%%%%%%%%%%%%%%%%%%%%%%%
\section{Proof of the Theorem}

\paragraph{An algebraic result.} 
Let $V$ be an arbitrary even lattice with an ordered basis $B=(e_1,\ldots,e_n)$ 
of roots and the corresponding Coxeter element $\tau$. For any root $a \in V$, 
we define two {\em Hilbert-Poincar\'e series} of $(V, a)$ by
\begin{align*}
P_{(V,a)}(t) &:= \sum_{k=0}^\infty \Big(1+\sum_{\ell=0}^{k-1} 
                  \langle a,\tau^\ell a\rangle\Big) t^k, \\
Q_{(V,a)}(t) &:= \sum_{k=0}^\infty \Big(1-\sum_{\ell=0}^{k} 
                  \langle a,\tau^{-\ell} a\rangle\Big) t^k.
\end{align*}
Define a form $(-,-)$ on $V$ by $(e_i,e_j)=-\langle e_i,e_j \rangle$ for $i<j$, 
$(e_i,e_j)=0$ for $i >j$, and $(e_i,e_i)=1$ for $1 \leq i,j \leq n$. Note that
 $\langle x,y \rangle = -(x,y)-(y,x)$ for all $x,y \in V$.
If $A$ denotes the matrix $((e_i,e_j))$ then, with respect to the basis $B$ of 
$V$, the mapping $\tau$ is given by the matrix $-A^{-1}A^t$. From this it follows 
that $(y,x)=-(x,\tau y)$. This implies, for a root $a\in V$,
\[
P_{(V,a)}(t) = \sum_{k=0}^\infty (a, \tau^k a)t^k, \qquad
Q_{(V,a)}(t) = \sum_{k=0}^\infty (a, \tau^{-k} a)t^k.
\]
These are the Hilbert-Poincar\'e series defined by Lenzing in \cite{L2}. By 
\cite[Lemma 18.1]{L2} these are always rational functions.

Now let $V_-$ be an arbitrary lattice with an ordered basis $B=(e_1,\ldots,e_n)$ 
of roots. Let $V_0:=V_-[u]$ and $V_+:=V_-[u,w]$. Let $\Delta_-, \Delta_0, \Delta_+$ 
be the characteristic polynomial of the Coxeter element of the lattice $V_-$, $V_0$,
$V_+$, respectively.

\begin{proposition} \label{prop:LP}
One has
\begin{align*}
{\rm (i)}  &&& Q_{(V_0,e_n)}(t) = \frac{\Delta_-(t)}{\Delta_0(t)},\\
{\rm (ii)} &&& P_{(V_0,e_n)}(t) + t =  \frac{\Delta_+(t)}{\Delta_0(t)}.
\end{align*}
\end{proposition}

\begin{proof}
The proof follows from elementary calculations.  Let $e_{n+1}:=e_n -u$. As
$u$ is in the radical of $V_0$, we have $P_{(V_0,e_n)}=P_{(V_0,e_{n+1})}$ and
$Q_{(V_0,e_n)}=Q_{(V_0,e_{n+1})}$. Now Part (i) of the proposition can be
derived from  the proof of \cite[Proposition~18.3]{L2}. Part (ii) follows from 
\cite[Corollary~18.2]{L2}.
\end{proof}

\paragraph{Geometrical model for the lattice.}
In this section, we will exhibit a geometrical model for the lattices
encountered in the previous section. Let $S$ be the smooth compactification
introduced before, and let $\bE$ be the tree of rational $(-2)$-curves.
Consider the K-group $K(S)=K(\Coh(S))$, which is generated by coherent sheaves
modulo relations coming from short exact sequences. Denote by $\Coh_\bE(S)$
the abelian subcategory of coherent sheaves on $S$ whose support is contained 
in $\bE$ and let $K_\bE(S)=K(\Coh_\bE(S))$ be its K-group. 

There is a natural bilinear pairing on K-groups, given by the Euler form, 
$\chi(A,B)=\sum_i(-1)^i\dim\Ext^i_S(A,B)$ for two coherent sheaves $A$ and $B$ 
on $S$. The sum is finite (actually only summands for $i=0,1,2$ occur) as $S$
is smooth, and each summand is finite as $S$ is projective. This form is not 
symmetric in general. Also note that 
 $\chi(\calO_S,B)=\chi(B)=\sum_i(-1)^i\dim H^i(S,B)$ 
is the usual Euler characteristic of a sheaf $B$ on $S$.

We denote by $N(S)$ the numerical K-group which is obtained from $K(S)$ by dividing
out the radical of the Euler form. As $S$ is a smooth and projective surface, $N(S)$
is a free abelian group, isomorphic to the algebraic part of cohomology, i.e.\
generated by the fundamental class, the N\'eron-Severi group and the class of a point. 
The identification of elements of the (numerical) K-group with cycles uses 
the Chern character $\chern: K(S)\to H^*(S,\ZZ)\otimes\QQ$. 
Note that only 2 can occur as denominator. However, we can just ignore this issue since
$\chern:K_\bE(S)\to H^*(S,\ZZ)$ is well-defined, as all components of $\bE$ have even
self-intersection and we have $\chern(\calO_D(l))=[D]+(l-D^2/2)\cdot pt$. Here, the
class of a point $pt$ is realised by the skyscraper sheaf $k(p)$ of a point, i.e.\
$pt=[\chern(k(p))]$. (Instead of cohomology, one could work equally well with the
Chow (intersection) ring as the target of the Chern character.)

Let $N_\bE(S)$ be the subgroup of $N(S)$ spanned by classes which have a 
representative with support in $\bE$. Thus, $N_\bE(S)$ is generated by the 
irreducible components of $\bE$ and the class of a point, $pt$. Lemma \ref{eulerform} 
shows that $pt$ is in the radical of $K_\bE(S)$. But it is not in the radical of
$K(S)$ by virtue of $\chi(\calO_S,k(p))=1$ and hence is a nonzero class in $N_\bE(S)$.

Using the notation of Figure~\ref{Fig1}, we introduce the following torsion sheaves 
with support contained in $\bE$:
\[ \F{j}{i} := \calO_{\E{j}{i}}(-1), \qquad
    F := \calO_E(-1), \qquad
    \widetilde F := \calO_E; \]
note $\chern(\calO_\F{j}{i}(-1))=[\E{j}{i}]$.
Their classes in $K_\bE(S)$ are denoted by
\[ \e{j}{i} := [\F{j}{i}] = [\calO_{\E{j}{i}}(-1)], \qquad
    e := [F] = [\calO_E(-1)], \qquad
    \widetilde e := [\widetilde F] := [\calO_{E}], \]
these form a basis. The class of the point $p\in\bE$ is given by 
$[k(p)]=\widetilde e-e$, in view of $0\to\calO_E(-1)\to\calO_E\to k(p)\to0$.

\begin{definition}
A coherent sheaf $G$ on $S$ is called \emph{spherical} if
\[ \Ext^l_S(G,G) = \left\{ \begin{array}{ll}
                                    \CC, & l=0 \text{ or } l=2 \\
                                    0    & \text{else}
                           \end{array} \right.
   \qquad \text{ and } \qquad
   G\otimes\omega_S\cong G .\]
\end{definition}
If $C$ is a complete, smooth, rational $(-2)$-curve on $S$, then the torsion sheaf
$\calO_C$ is spherical: Firstly, we have $\Hom_S(\calO_C,\calO_C)=\CC$ since $C$ is 
connected and $\Ext^1_S(\calO_C,\calO_C)=0$ as $C$ is rigid. The adjunction formula 
shows 
$-2=2g_C-2=\deg(K_C)=\deg((K_S\otimes\calO_S(C))|_C)=\deg(K_S|_C)-2$,
hence $\deg(K_S|_C)=0$ and $\omega_S|_C\cong\calO_C$. Now Serre duality shows that
$\dim(\Ext^2_S(\calO_C,\calO_C))=\dim(\Hom_S(\calO_C,\calO_C))=1$. By a similar
computation, the sheaves $\calO_C(k)$ are spherical as well, for any $k\in\ZZ$.

In particular, the sheaves $\F{j}{i}$ and $F$ and $\widetilde F$ are all spherical.
At this point we note that the Euler form on $K_\bE(S)$ and $N_\bE(S)$ is symmetric 
--- this is a direct consequence of the second sphericality condition and Serre 
duality. 

\bigskip

We define the lattices
\begin{align*}
 V_0 &:= N_\bE(S)\subset N(S), && \langle-,-\rangle := -\chi(-,-), \\
 V_- &:= N_\bE(S)\cap[\calO_S]^\perp.
\end{align*}
Letting $u:=pt=\widetilde e-e$, the following lemma shows $V_0=V_-[u]$ and that
$B_-=(\e{1}{1},\dots,\e{\alpha_r-1}{r},e)$ and $B_0=(B_-,e-u)$ are bases of roots.
\begin{lemma} \label{eulerform}
Let $C$, $C'$, $C''$ be smooth, rational $(-2)$-curves in a projective
surface $S$, with $C$ and $C'$ intersecting transversally and $C$ disjoint
to $C''$. Also, let $p\in C$ be a point. Then, for any $l\in\ZZ$,
\begin{align*}
(0) &&  0 &= \langle k(p),k(p)\rangle. &&&
(1) &&  1 &= \langle\calO_C(l),\calO_{C'}(l)\rangle. \\
(2) &&  0 &= \langle k(p), \calO_C(l)\rangle. &&&
(3) && -2 &= \langle\calO_C(l),\calO_C(l)\rangle. \\
(4) &&  0 &= \langle\calO_C(l),\calO_{C''}(l)\rangle. &&&
(5) && -2 &= \langle\calO_C(l),\calO_C(l+1)\rangle. \\
(6) &&  0 &= \langle\calO_S,\calO_C(-1)\rangle.
\end{align*}
\end{lemma}

\begin{proof}
(6) follows from $H^*(S,\calO_C(-1))=H^*(\PP^1,\calO_{\PP^1}(-1))=0$.
Applied to the lattices above, this implies $V_-\subset V_0$.

(0) The class $u=[k(p)]$ is isotropic in view of 
$\Hom_S(k(p),k(p))=\CC$, $\Ext^1_S(k(p),k(p))=\CC^2$ and
$\Ext^2_S(k(p),k(p))=\CC$. (2) The exact sequence 
$0\to\calO_S(-C)\to\calO_S\to\calO_C\to0$ shows
 $\chi(\calO_C,k(p)) = \chi(\calO_S,k(p))-\chi(\calO_S(-C),k(p)) = 1-1$.
These two properties show that $V_0=V_-\oplus\ZZ u$ is an orthogonal 
decomposition with isotropic $u$, i.e.\ $V_0=V_-[u]$.
 
(3) follows directly from the first condition of sphericality.

For (1), we have $\Hom_S(\calO_C,\calO_{C'})=0$. Serre duality and $\omega_S|_C\cong\calO_C$
imply $\Ext^2_S(\calO_C,\calO_{C'})=0$. Then, 
$\dim\Ext^1_S(\calO_C,\calO_{C'})=-\chi(\calO_C,\calO_{C'})=1$. The
last equation uses the ideal sheaf sequence for $C$ again, together with
$\chi(\calO_S(-C),\calO_{C'})=\chi(\calO(C)|_{C'})=\chi(\calO_{\PP^1}(1))=2$
and $\chi(\calO_S,\calO_{C'})=1$. The same computation works for arbitrary $l$.

(4) This is obvious as the supports are disjoint.

(5) follows directly from (3) and the fact that $[k(p)]$ is isotropic.
\end{proof}

% The curves $E_1^{(j)},\dots,E_{\alpha_j-1}^{(j)}$ form an $A_{\alpha_j-1}$-chain and, 
% by the lemma, the associated $\langle-,-\rangle$-matrix is
% \[ \calB^{(j)} := 
% \left(\begin{array}{rrrrr}
% -2&  1 & 0 &\dots& 0 \\
%  1& -2 & 1 &     & \vdots \\
%  0&  1 & \ddots & \ddots& 0 \\
%  \vdots & &\ddots&  -2 & 1 \\ 
%  0& \dots  & 0 &   1 & -2
% \end{array} \right) \]

% With respect to the ordered basis $(e_1^{(1)},\dots,e_{\alpha_r-1}^{(r)},e,\widetilde e)$
% of $V_0=N_\bE(S)$, we obtain the following $\langle-,-\rangle$-matrix
% \[ \calB := \left(\begin{array}{cccccc}
% \calB^{(1)} & 0 & \dots & 0 & {}^tb^{(1)} & {}^tb^{(1)} \\
% 0 & \calB^{(2)} & \dots & 0 & {}^tb^{(2)} & {}^tb^{(2)} \\
% \vdots & \vdots& \ddots& \vdots & \vdots     & \vdots \\
% 0      & 0 &   & \calB^{(r)}& {}^tb^{(r)} & {}^tb^{(r)} \\
% b^{(1)} & b^{(2)}& \dots & b^{(r)} & -2 & -2 \\         
% b^{(1)} & b^{(2)}& \dots & b^{(r)} & -2 & -2 \\         
% \end{array} \right) \text{ with }
% b^{(j)}=(0,\dots,0,1)\in\ZZ^{\alpha_j-1}
% \]
%
% First of all, this is just the matrix associated to the graph in Figure~\ref{Fig2}. 
% Next, the matrix $\calA$ obtained from $\calB$ by replacing the diagonal entries 
% by $-1$ and by zeroing the upper right part satisfies $\calB=\calA+{}^t\calA$. 
% This implies that there is a bilinear form $(-,-)$ on $N_\bE(S)$ such that 
% $\langle x,y\rangle=-(x,y)-(y,x)$. The matrix of $(-,-)$ is $-\calA$ and it is 
% triangular with all diagonal entries being 1. Hence we are free to use the results 
% from Section~\ref{SectionBila}.

\paragraph{The Fuchsian case.}
In this section, we explain how the geometric picture of the last section can be 
extended to certain Fuchsian singularities.

Let $(X,x)$ be a normal surface singularity with a good $\CC^*$-action. Denote by 
$p:\calX\to\calS$ the semiuniversal deformation. The $\CC^*$-action extends, at
least formally, to this deformation. According to the weights, we get subspaces 
$\calS^+$, $\calS^0$ and $\calS^-$ of $\calS$. These subspaces are unique up to
non-unique isomorphism. By pullback, we obtain a family $p^-:\calX^-\to\calS^-$. We
denote by $p^c:\calY\to\calS^-$ the simultaneous resolution of the singularities 
at infinity in all fibres after compactification; cf.\ e.g. \cite{Pinkham78}. 
In other words, all fibres of $p^c$ are projective surfaces having at most the 
vertices as singularities.

Now suppose that $(X,x)$ is Fuchsian with $g=0$. Then the special fibre of $p^c$ 
contains the configuration $\bE$. Actually, the configuration deforms, so that $\bE$ 
appears in every fibre. Next, assume that $(X,x)$ is \emph{negatively smoothable}, 
i.e.\ that the generic fibre of $p^c$ is smooth. By \cite[Proposition 6.13]{Pinkham78}, this 
implies that the generic fibre is a smooth K3 surface which we denote by $\genfib$. 
For example, if $(X,x)$ is an isolated hypersurface or complete intersection 
singularity, then it is negatively smoothable. On the other hand, since the rank 
of the N\'{e}ron-Severi group of a K3 surface is 20, a necessary condition for 
negative smoothability is $\sum_i\alpha_i\leq19+r$. 

In this situation, we will use $\genfib$ instead of $S$. Both are smooth projective
surfaces containing $\bE$, but $\genfib$ has the advantage of being a K3 surface 
which means that the structure sheaf $\calO_\genfib$ is spherical. Equipping the 
K-groups $K(\genfib)$ and $K_\bE(\genfib)$ with the Euler form, we introduce the lattice
\[
  V_+ := N_\bE(\genfib)\oplus\ZZ[\calO_\genfib]\subset N(\genfib), \qquad
            \text{with } w=u-[\calO_\genfib].
\]
Here, $-w=[\calO_\genfib]-[k(p)]=[\calI_p]$ is the class of the ideal sheaf corresponding 
to a point $p$ of $\bE$. As $V_-$ and $V_0$ did before, $V_+$ has a basis of (classes of) 
spherical sheaves. Serre duality, coupled with $\omega_\genfib\cong\calO_\genfib$,
shows that the Euler form is symmetric. Also note that the lattices $V_-$ and $V_0$ 
can be considered as sublattices of $V_+$. The vectors $u$ and $w$ do indeed form a 
hyperbolic plane which is orthogonal to $V_-$ by virtue of
\begin{align*}
 (7) &&  0 &= \langle\calO_\genfib,\calO_C(-1)\rangle. &&&
 (8) && -2 &= \langle\calO_\genfib,\calO_\genfib\rangle.
\end{align*}

\paragraph{Computing the Poincar\'e series.}
Theorem \ref{thm:main} immediately follows from the following fact, coupled with 
Proposition \ref{prop:LP}.
\begin{proposition}
The series 
$Q_{(V_0, e)}(t)$ and  $P_{(V_0, e)}(t) + t$ are the Poincar\'e series of a Kleinian 
singularity  and a Fuchsian singularity of genus 0, respectively.
\end{proposition}

\begin{proof}
The element $u=\widetilde{e} -e$ spans the radical ${\rm rad}\, V_0$ of the lattice $V_0$. 
Let $\overline{V}_0= V_0/{\rm rad}\, V_0$. For an automorphism $\sigma$ of $V_0$, we
denote by the same letter the induced automorphism
 $\sigma : \overline{V}_0 \to \overline{V}_0$.
By more abuse of notation, we will denote elements of $V_0$ and their classes in 
$\overline{V}_0$ by the same letter. In order to compute the series $P_{(V_0, e)}(t)$ and 
$Q_{(V_0, e)}(t)$, it suffices to consider the automorphism $\tau_0$ of $\overline{V}_0$. 
Now one can easily see that, on the quotient space, $s_e s_{e-u}={\rm id}_{\overline{V}_0}$. 
Therefore we have, again on the quotient $\overline{V}_0$,
\[ \tau_0 = \tau_1 \cdots \tau_r \mbox{ where } 
   \tau_i:= s_{\e{1}{i}} \cdots s_{\e{\alpha_i-1}{i}} \mbox{ for } \ i=1, \ldots , r. \]
For $\tau_i$ we have
 $\tau_i(e)=e+ \sum_{j=1}^{\alpha_i-1} \be{j}{i}$, 
 $\tau_i(\be{j}{i})=\be{j+1}{i}$ 
for $j=1, \ldots , \alpha_i-2$, 
$\tau_i(\be{\alpha_i-1}{i})= - \sum_{j=1}^{\alpha_i-1} \be{j}{i}$, and $\tau_i$ 
is the identity on all other basis elements of $\overline{V}_0$. Therefore (if $\alpha_i > 2$)
\[ \tau_i(e)=e+ \sum_{j=1}^{\alpha_i-1} \be{j}{i}, \quad
   \tau_i^2(e)= e+ \sum_{j=2}^{\alpha_i-1} \be{j}{i}, \quad 
   \ldots, \quad
   \tau_i^{\alpha_i}(e)=e. \]
From this it follows that for $k>0$
\begin{equation}  \label{FuchsDiv}
   1 + \langle e, \sum_{\ell=0}^{k-1} \tau_0^\ell e \rangle 
=  1 + \langle e, ke 
     + \sum_{i=1}^r \left[ k \frac{\alpha_i - 1}{\alpha_i} \right]  \e{\alpha_i -1}{i} \rangle
=  1 + \deg D_{\rm Fuchs}^{(k)}.
\end{equation}
By Riemann-Roch, since $r \geq 3$  in the Fuchsian case
\[ \dim L(D_{\rm Fuchs}^{(k)}) = 1+ \deg D_{\rm Fuchs}^{(k)} \mbox{ for } k \neq 1, \]
and $\dim L(D_{\rm Fuchs}^{(1)}) = \dim L(D_0)=g=0$.

Similarly we have in the Kleinian case,
 $\tau_0^{-1} = \tau_r^{-1} \cdots \tau_1^{-1}$, 
 $\tau_i^{-1}(e) = \overline{e} + \be{\alpha_i-1}{i}$, 
 $\tau_i^{-1}(\be{1}{i})= - \sum_{j=1}^{\alpha_i-1} \be{j}{i}$, and 
 $\tau_i^{-1}(\be{j}{i})=\be{j-1}{i}$ for $j=2, \ldots , \alpha_i-1$. 
Therefore (if $\alpha_i >2$)
\[ \tau_i^{-1}(e) = e + \be{\alpha_i-1}{i}, \quad
   \ldots, \quad
   \tau_i^{1-\alpha_i}(e)=e+\sum_{j=1}^{\alpha_i-1}\e{j}{i}, \quad
   \tau_i^{-\alpha_i}(e)=e.\] 
From this it follows that for $k > 0$
\begin{equation}  \label{KleinDiv}
   1 - \langle e, \sum_{\ell=1}^k \tau_0^{-\ell} e \rangle 
=  1 - \langle e, ke 
     + \sum_{i=1}^r  \left(k - \left[ \frac{k}{\alpha_i} \right] \right) \e{\alpha_i -1}{i} \rangle 
=  1 + \deg D_{\rm Klein}^{(k)}.
\end{equation}
Finally, since $r \leq 3$ in the Kleinian case, Riemann-Roch gives
\[ \dim L(D_{\rm Klein}^{(k)}) = 1+ \deg D_{\rm Klein}^{(k)}. \]
We want to point out that the righthand equality in (\ref{KleinDiv}) is not just 
a numerical coincidence. Rather, it follows from the fact that the restriction of 
the divisor  $kE + \sum_i(k-[ k/\alpha_i])\E{\alpha_i -1}{i}$
to $E$ is precisely $D_{\rm Klein}^{(k)}$; similarly for the Fuchsian case 
and equation (\ref{FuchsDiv}).
\end{proof}

%%%%%%%%%%%%%%%%%%%%%%%%%%%%%%%%%%%%%
\section{Categorical interpretation}

\paragraph{The Coxeter functor.}
Denote by $D^b(S)$ the bounded derived category of coherent sheaves on $S$, 
and by $D^b_\bE(S)$ the subcategory of complexes with support (of all
homology sheaves) contained in $\bE$. These are triangulated categories 
which should be seen as geometrical instances of homological invariants
like (numerical) K-groups. A first hint at this is the fact that to any 
spherical sheaf $G$ one can associate a \emph{spherical twist} 
$\TTT_G:D^b(S)\to D^b(S)$ which is an autoequivalence. It is defined by 
distinguished triangles 
$\Hom^\bullet_S(G,A)\otimes_\CC G\to A\to \TTT_G(A)$ for $A\in D^b(S)$. 
See \cite[\S5.3]{Huybrechts} or \cite{ST} for details. 
For a spherical sheaf $G$ with support in $\bE$, the functor $\TTT_G$ 
induces an autoequivalence of $D^b_\bE(S)$.

To any triangulated category one can associate its K-group. In the case 
of $D^b(S)$ and $D^b_\bE(S)$, these are $K(S)$ and $K_\bE(S)$, respectively. 
The triangulated category  $D^b_\bE(S)$ is 2-Calabi-Yau, i.e.\ has the even 
shift $[2]$ as Serre functor, and this explains why the Euler form on $K_\bE(S)$
is symmetric. (The Serre functor of $D^b(S)$ is $-\otimes\omega_S[2]$ 
and thus is in general not a shift.)

Any autoequivalence $\varphi:D^b(S)\to D^b(S)$ induces an automorphism 
$\varphi^K:K(S)\to K(S)$ of the K-group, which is an isometry for the 
Euler pairing. In case $\varphi=\TTT_G$ is a spherical twist, we have 
$\chi(G,G)=2$ as $S$ is a surface. The defining triangles for $\TTT_G$ 
show that $\TTT_G^K$ is just the reflection along $[G]\in K(S)$.

This allows us to lift lattices and Coxeter elements to the categorical 
level (the third lines about $\calD_+$ and $\varphi_+$ only apply in the
case of negatively smoothable Fuchsian singularities): 
\begin{align*}
\calD_- &:= D^b_\bE(S)\cap\calO_S^\perp, &&& V_- &= N(\calD_-),\\
\calD_0 \, &:= D^b_\bE(S),                 &&& V_0 \, &= N(\calD_0),\\
\calD_+ &:= \langle D^b_\bE(\genfib),\calO_\genfib\rangle, &&& V_+ &= N(\calD_+),
\intertext{and}
\varphi_- &:= \TTT_{\calO_{E_1^{1}}(-1)} \cdots
              \TTT_{\calO_{E_{\alpha_r-1}^{r}}(-1)} \TTT_{\calO_E(-1)},
                                      &&& \tau_- &= \varphi_-^K, \\
\varphi_0 \, &:= \TTT_{\calO_{E_1^{1}}(-1)} \cdots
              \TTT_{\calO_{E_{\alpha_r-1}^{r}}(-1)} \TTT_{\calO_E(-1)} \TTT_{\calO_E},
                                      &&& \tau_0 \; &= \varphi_0^K, \\
\varphi_+ &:= \TTT_{\calO_{E_1^{1}}(-1)} \cdots
              \TTT_{\calO_{E_{\alpha_r-1}^{r}}(-1)} 
              \TTT_{\calO_E(-1)} \TTT_{\calO_E} \TTT_{\calO_\genfib},
                                      &&& \tau_+ &= \varphi_+^K. 
\end{align*}
At the top, $\calO_S^\perp:=\{A\in D^b(S)\:|\:\Hom^*(\calO_S,A)=0\}$. 
Hence $\calD_-$ and $\calD_0$ are triangulated subcategories of $D^b(S)$. 
$\calD_+$ is, by definition, the smallest triangulated subcategory of $D^b(\genfib)$
containing $D^b_\bE(\genfib)$ and $\calO_\genfib$.

Turning to the K-groups, it is clear that $K(\calD_0)=K_\bE(S)$ and 
$K(\calD_-)=[\calO_S]^\perp$ and $K(\calD_+)=K_\bE(\genfib)\oplus\ZZ[\calO_\genfib]$.
The numerical K-groups are all defined using the radical of the K-group of the
surface. For example, $N(\calD_-)$ is the image of the composition
 $K(\calD_-)\hookrightarrow K(D^b(S))=K(S)\to N(S)$; similar for the other
lattices which are all equipped with the negative of the Euler form.

As a side remark, we want to point out that the three categories can also be 
described by $\calD_0=\langle B_0\rangle$, $\calD_-=\langle B_-\rangle$, and 
$\calD_+=\langle B_+\rangle$, i.e.\ they coincide with the smallest triangulated
subcategories containing the basis sheaves. For this, one uses that $\bE$ 
consists of rigid projective lines.\footnote{
For any point $p\in E$, we have $k(p)\in\langle B_0\rangle$ as the cokernel of a morphism
$\calO_E(-1)\to\calO_E$. This implies $k(p')\in\langle B_0\rangle$ for all $p'\in\bE$ since
$\bE$ is a tree of rational curves. Hence, if $A\in\Coh(S)$ is supported on one component 
of $\bE$, then $A\in\langle B_0\rangle$. Note, that there are no non-trivial self-extensions 
with one-dimensional support as the irreducible components are rigid. Finally, an arbitrary 
sheaf $A\in\Coh(S)$ with support in $\bE$ can be split into pieces with support on different 
components. For this, use the cones of $A\to i_*{\rm L}i^*A$ for the inclusion 
$i:E\hookrightarrow S$. Then, $\Coh_\bE(S)\subset\langle B_0\rangle$ immediately implies 
$D^b_\bE(S)=\langle B_0\rangle$. The reasoning for the other categories proceeds along the same 
lines. It is worth mentioning that $D^b(\bE)$ does not have finite generators --- it is 
representation-wild; see \cite{DG}.
}
Each of the three triangulated categories is generated by spherical objects
and the autoequivalences $\varphi_-$, $\varphi_0$ and $\varphi_+$ are given
by a concatenation of spherical twist functors. Their traces on the K-groups
are precisely the Coxeter elements we used before.

\paragraph{Another approach by Kajiura, Saito and Takahashi.}
They study the triangulated category $\Dsgr(A)$ of graded singularities of the algebra 
$A$ in the hypersurface case. It is defined as the quotient of the bounded derived 
category of finitely generated graded $A$-modules by the derived category of projective 
modules. Well known other descriptions for $\Dsgr(A)$ use either graded maximal 
Cohen-Macaulay modules or, in the hypersurface case, matrix factorisations.

In \cite{KST1} and \cite{KST2}, the authors prove, for special cases of $A$, the existence 
of full, strongly exceptional sequences in $\Dsgr(A)$. These yield root bases in the 
associated Grothendieck groups (with respect to the symmetrised Euler form) and a corresponding
Coxeter element. The latter lifts to an autoequivalence of $\Dsgr(A)$ given by the
Auslander-Reiten translation functor.

\cite{KST1} addresses regular weight systems with $\varepsilon=1$, which correspond 
precisely to the Kleinian singularities (using the notation from Section \ref{KandF}, 
one has $\varepsilon=-R$). \cite{KST2} deals with certain regular weight systems with 
$\varepsilon=-1$, which are the hypersurface singularities of Fuchsian type with $g=0$
and either $r=3$ (the 14 unimodal singularities) or $r=4$ (six bimodal singularities) or
$r=5$ (two singularities). However, the virtue of the approach presented here seems to be
that the Coxeter element lifts to an autoequivalence of the derived category of the ambient
surfaces; also, we don't need to impose any restriction on the Fuchsian singularities (apart
from negative smoothability for the categorical interpretation).

%%%%%%%%%%%%%%%%%%%%%%%%%%%%%%%%%%%%%

\bigskip
\noindent Leibniz Universit\"{a}t Hannover, Institut f\"{u}r Algebraische Geometrie,\\
Postfach 6009, D-30060 Hannover, Germany \\
E-mail: ebeling@math.uni-hannover.de\\
E-mail: ploog@math.uni-hannover.de

\end{document}